\documentclass[smallextended]{svjour3}       
\usepackage{amsmath,amssymb,hyperref,color}

\smartqed  
\usepackage{graphicx}
\journalname{Journal} 
\begin{document}

\title{Spatial Search on Johnson Graphs by Continuous-Time Quantum Walk
}


\author{
        Hajime Tanaka$^1$           \and
        Mohamed Sabri$^1$         \and
        Renato Portugal$^2$ 
}

\authorrunning{Hajime Tanaka \and Mohamed Sabri \and  Renato Portugal} 

\institute{$^1$H. Tanaka and M. Sabri      \at
              Research Center for Pure and Applied Mathematics\\
              Graduate School of Information Sciences\\
              Tohoku University\\
              Sendai 980-8579, Japan \\
              \email{htanaka@tohoku.ac.jp, sabrimath@dc.tohoku.ac.jp}           
           \and
           $^2$R. Portugal              \at
           National Laboratory of Scientific Computing (LNCC) \\
           Petr\'opolis, RJ, 25651-075, Brazil\\
           \email{portugal@lncc.br}
}


\maketitle

\begin{abstract}
Spatial search on graphs is one of the most important algorithmic applications of quantum walks. To show that a quantum-walk-based search is more efficient than a random-walk-based search is a difficult problem, which has been addressed in several ways. Usually, graph symmetries aid in the calculation of the algorithm's computational complexity, and Johnson graphs are an interesting class regarding symmetries because they are regular, Hamilton-connected, vertex- and distance-transitive. In this work, we show that spatial search on Johnson graphs by continuous-time quantum walk achieves the Grover lower bound $\pi\sqrt{N}/2$ with success probability $1$ asymptotically for every fixed diameter, where $N$ is the number of vertices. The proof is mathematically rigorous and can be used for other graph classes.

\end{abstract}

\section{Introduction}

The continuous-time quantum walk (CTQW) was introduced by Farhi and Gutmann~\cite{FG98} as a quantum analogue of the continuous-time Markov process with the aim of building faster quantum algorithms based on decision trees. An example of such algorithm was provided for the problem of evaluating NAND-based Boolean formulas: if the NAND tree has $N$ leaves, the best classical algorithm runs in time $\Omega(N^{0.753})$, and the quantum-walk-based algorithm runs in time $O(N^{0.5})$~\cite{FGG_08}. Childs \textit{et al.}~\cite{CCDFGS_03} also presented an interesting example of an exponential algorithmic speedup by using the CTQW to rapidly traverse a graph.

Besides algorithmic applications, many other important results were obtained. We highlight some of them.
Konno~\cite{Kon_05,Kon_06} proved a weak limit theorem for the CTQW on the line and trees, and showed that there is a striking contrast to the central limit theorem of symmetric classical random walks.
M\"ulken and Blumen~\cite{MB_11} reviewed applications of CTQWs to transport in various systems, and recently Razzoli \textit{et al.}~\cite{RPB_21} analytically determined subspaces of states having maximum transport efficiency for many graph topologies. Benedetti \textit{et al.}~\cite{BRP_19} described CTQW on dynamical percolation graphs with the goal of analyzing the effect of noise produced by randomly adding or removing graph edges.
Delvecchio \textit{et al.}~\cite{DPW_20} analytically investigated the analogy between the CTQW in one dimension and the evolution of the quantum kicked rotor at quantum resonance conditions.

Benioff~\cite{Ben02} originally proposed the quantum spatial search problem, which aims to find a marked location using a \textit{quantum robot} wandering aimlessly on a graph endowed with the skill of checking whether a node is marked or not. The continuous-time model proved useful to solve this problem on many graphs. Childs and Goldstone~\cite{CG_04} presented a framework using a Hamiltonian that encodes the adjacency matrix and the information about the location of the marked vertex. This framework has been used systematically in many papers; we highlight some of them. Agliari \textit{et al.}~\cite{ABM_10} studied the spatial search on fractal structures and showed how the transition from the ground state of the Hamiltonian to a state close to the marked state is accomplished by a CTQW. 
Philipp \textit{et al.}~\cite{PTB_16} analyzed spatial search on balanced trees, and showed that the efficiency depends on whether the marked vertex is close to the root or the leaves of the graph.
Osada \textit{et al.}~\cite{OCOSMN_20} explored the spatial search on scale-free networks and found that the efficiency is determined by the global structure around the marked vertex.
Recently, new experimental implementations of quantum walk search were proposed in~\cite{DGPSW_20,DGGCWS_19}.

The Johnson graph $J(n,k)$ has many interesting properties that help to establish the efficiency of quantum-walk-based search. Johnson graphs are regular, vertex- and distance-transitive. The eigenvalues and eigenvectors of their adjacency matrices are known and easy to handle. In fact, there are some results in literature analyzing quantum walk search on a subclass of Johnson graphs. Wong~\cite{Wong2016JPA} showed that $J(n,3)$ supports fast spatial search by continuous-time quantum walk, and Xue \textit{et al.}~\cite{XRL_19} also showed that $J(n,3)$ supports fast spatial search using the scattering quantum walk, and they discussed the general case when $k$ is arbitrary. The scattering quantum walk is a discrete-time model equivalent to the coined model~\cite{AL_09}.

In this work, we use the continuous-time model to analyze quantum-walk-based search algorithms on the Johnson graph $J(n,k)$ for arbitrary $n$ and $k$. We show that the optimal running time is $\pi\sqrt{N}/2$ and the success probability is $1$ asymptotically on $J(n,k)$ for every fixed $k$, where $N=\binom{n}{k}$ is the number of vertices. We employ a fully rigorous analytic method, which can be borrowed for other graph classes.

The outline of this paper is as follows. In Sec.~\ref{sec:frame}, we lay out the mathematical framework that we use to prove the computational complexity of the search algorithm. In Sec.~\ref{sec:search}, we show that the quantum-walk-based search achieves the lower bound on Johnson graphs. In Sec.~\ref{sec:conc}, we draw our conclusions.

\section{Mathematical framework}\label{sec:frame}

Let $J(n,k)$ be the Johnson graph, where the vertex set $V=V(J(n,k))$ is the set of $k$-subsets of $[n]=\{1,2,\dots,n\}$, and two vertices $v,v'\in V$ are adjacent if and only if $|v\cap v'|=k-1$.
Note that $J(n,1)$ is the complete graph $K_n$, and that $J(n,2)$ is the triangular graph $T_n$, which is strongly regular. We associate $J(n,k)$ with a Hilbert space spanned by $\mathcal{H}=\{|v\rangle : |v\rangle\in V\}$, as is usually done in the definition of the continuous-time quantum walk~\cite{FG98}.

Let $w\in V$ be the marked vertex.
We consider the Hamiltonian of the form~\cite{CG_04}
\begin{equation*}
	H=-\gamma A- | w\rangle \langle w |,
\end{equation*}
where $A$ denotes the adjacency operator of $J(n,k)$, and $\gamma$ is a real and positive parameter.
The algorithm starts in the initial state $|\psi(0)\rangle$, which is the uniform superposition of the computational basis
\begin{equation*}
	|{{s}}\rangle:=\frac{1}{\sqrt{N}}\sum_{v\in V}|v\rangle,
\end{equation*}
where $N=\binom{n}{k}$ is the number of vertices of $J(n,k)$, and the notation ``:='' is used to stress that we are defining a new symbol.
The quantum state at time $t$ is therefore given by
\begin{equation*}
	|\psi(t)\rangle=\mathrm{e}^{-\mathrm{i}Ht}|{{s}}\rangle.
\end{equation*}

\subsection{Invariant subspace}
From now on, we fix $k$. We
will always assume that $n\geqslant 2k$, so that $J(n,k)$ has diameter $k$.
Clearly, this assumption does not affect the asymptotic analysis of the search algorithm when $n\rightarrow\infty$.
Consider the following subsets of $V$
\begin{equation*}
	\nu_{\ell}=\{v\in V:|v\cap w|=k-\ell\},
\end{equation*}
where $0\leqslant \ell\leqslant k$. Note that $\nu_0=\{w\}$, and that $|\nu_{\ell}|=\binom{k}{\ell}\binom{n-k}{\ell}$.
Instead of  using the Hilbert space $\operatorname{span}\{|v\rangle : v\in V\}$, we will work with the invariant subspace $\mathcal{H}_{\mathrm{inv}}=\operatorname{span}\{|\nu_{\ell}\rangle : 0\leqslant \ell\leqslant k\}$, where
\begin{equation*}
	|\nu_{\ell}\rangle=\frac{1}{\sqrt{|\nu_{\ell}|}} \sum_{v\in\nu_{\ell}} |v\rangle,
\end{equation*}
for $0\leqslant \ell\leqslant k$.

Below we collect the necessary information regarding Johnson graphs~\cite{BCN1989B,DKT_18}.
The adjacency operator $A$ has exactly $k+1$ distinct eigenvalues $\lambda_0>\lambda_1>\dots>\lambda_k$, where
\begin{equation}\label{eigenvalues}
	\lambda_{\ell}=(k-\ell)(n-k-\ell)-\ell,
\end{equation}
and the multiplicity of $\lambda_{\ell}$ equals $\binom{n}{\ell}-\binom{n}{\ell-1}$ (with the understanding that $\binom{n}{-1}=0$).
For $0\leqslant \ell\leqslant k$, let $P_{\ell}$ denote the projector onto the eigenspace of $A$ in $\operatorname{span}\{|v\rangle : v\in V\}$ for the eigenvalue $\lambda_{\ell}$.
It is known that
\begin{equation}\label{Pw}
	\| P_{\ell} | w \rangle \|^2 = \langle w | P_{\ell} | w \rangle= \frac{\binom{n}{\ell}-\binom{n}{\ell-1}}{\binom{n}{k}}=\frac{k!(n-k)!(n-2\ell+1)}{\ell!(n-\ell+1)!}.
\end{equation}
In particular, the vectors $P_{\ell} | w \rangle$ are nonzero, and hence $\mathcal{H}_{\mathrm{inv}}$ has another orthonormal basis $\{|\lambda_{\ell}\rangle : 0\leqslant \ell\leqslant k\}$ consisting of eigenvectors of $A$, where
\begin{equation*}
	|\lambda_{\ell}\rangle = \frac{1}{\| P_{\ell} | w \rangle \|} P_{\ell} | w \rangle.
\end{equation*}
We will analyze the search algorithm using this basis.
We note that this is different from the basis used by Wong~\cite{Wong2016JPA}, who analyzed the case $k=3$.

Observe that $P_0=|{{s}}\rangle\langle{{s}}|$, from which it follows that
\begin{equation*}
	|\lambda_0\rangle = |{{s}}\rangle.
\end{equation*}
The matrix representation of $A$ in this basis is
\begin{equation}\label{A}
	A=\operatorname{diag}(\lambda_0,\lambda_1,\dots,\lambda_k).
\end{equation}
On the other hand, since
\begin{eqnarray*}
	\big(|w\rangle\langle w|\big)|\lambda_{\ell}\rangle &=& \| P_{\ell}|w\rangle\| |w\rangle = \| P_{\ell}|w\rangle\|\sum_{\ell'=0}^kP_{\ell'}|w\rangle\\
	&=& \| P_{\ell}|w\rangle\|\sum_{\ell'=0}^k \| P_{\ell'}|w\rangle \| |\lambda_{\ell'}\rangle,
\end{eqnarray*}
the matrix representation of $|w\rangle\langle w|$ in this basis is
\begin{equation}\label{w}
	|w\rangle\langle w|=\big( \| P_{\ell}|w\rangle\| \| P_{\ell'}|w\rangle\| \big)_{\ell,\ell'=0}^k.
\end{equation}

\subsection{Degenerate perturbation theory}

We apply the degenerate perturbation theory to the analysis of the search algorithm~\cite{JMW_14,Wong2016JPA}.
In perturbation theory, eigenvalues and eigenvectors of an analytic square matrix function are \emph{assumed} to be again analytic functions. However, this seems to be a subtle problem, particularly for degenerate eigenvalues, i.e., eigenvalues with multiplicity at least two. To make our discussions \emph{fully} rigorous, we proceed as follows.

Set
\begin{equation*}
	\epsilon:=\frac{1}{\sqrt{n}},
\end{equation*}
and view the parameter $\gamma$ as a function of $\epsilon$ (where we recall $k$ is fixed).
Let
\begin{equation*}
	\eta:=\frac{1}{\gamma n}=\frac{\epsilon^2}{\gamma},
\end{equation*}
so that
\begin{equation*}
	-\eta H=\epsilon^2 A+\eta |w\rangle\langle w|.
\end{equation*}
From \eqref{eigenvalues} it follows that
\begin{equation*}
	r_{\ell}(\epsilon):=\epsilon^2\lambda_{\ell}=(k-\ell)(1-(k+\ell)\epsilon^2)-\ell\epsilon^2,
\end{equation*}
where $0\leqslant \ell\leqslant k$.
Likewise, from \eqref{Pw} it follows that
\begin{equation*}
	p_{\ell}(\epsilon) := \| P_{\ell}|w\rangle \| = \epsilon^{k-\ell}\sqrt{\frac{ k! (1-(2\ell-1)\epsilon^2) }{ \ell! (1-(\ell-1)\epsilon^2)\cdots (1-(k-1)\epsilon^2) }}.
\end{equation*}
By these comments and \eqref{A} and \eqref{w}, the matrix representation of $-\eta H$ in the above basis is
\begin{equation*}
	-\eta H=\operatorname{diag}(r_0(\epsilon),r_1(\epsilon),\dots,r_k(\epsilon)) +\eta \big(p_{\ell}(\epsilon)p_{\ell'}(\epsilon)\big)_{\ell,\ell'=0}^k.
\end{equation*}

Our aim is to invoke the implicit function theorem for complex analytic functions~\cite{Krantz1992B}.
To this end, we extend for the moment the range of $\epsilon$ to complex numbers with $|\epsilon|^2<(2k-1)^{-1}$, so that the functions $r_{\ell}(\epsilon)$ and $p_{\ell}(\epsilon)$ are all analytic.
We now fix $a\in\mathbb{C}\setminus\{0\}$, and consider the following $k+3$ analytic functions of $k+4$ variables:
\begin{align*}
	f_0(\epsilon,\eta,\xi_0,\xi_1,\dots,\xi_k,\lambda) &= r_0(\epsilon)\xi_0+\eta\, p_0(\epsilon)\sum_{\ell=0}^k p_{\ell}(\epsilon) \xi_{\ell} -\lambda \xi_0, \\
	& \ \, \vdots \\
	f_k(\epsilon,\eta,\xi_0,\xi_1,\dots,\xi_k,\lambda) &= r_k(\epsilon)\xi_k+\eta\,p_k(\epsilon)\sum_{\ell=0}^k p_{\ell}(\epsilon) \xi_{\ell} -\lambda \xi_k, \\
	f_{k+1}(\epsilon,\eta,\xi_0,\xi_1,\dots,\xi_k,\lambda) &= \xi_0 -a, \\
	f_{k+2}(\epsilon,\eta,\xi_0,\xi_1,\dots,\xi_k,\lambda) &= \sum_{\ell=0}^k p_{\ell}(\epsilon) \xi_{\ell} -1.
\end{align*}
Note that $f_j(\epsilon,\eta,\xi_0,\xi_1,\dots,\xi_k,\lambda)=0$ for $0\leqslant j\leqslant k$ if and only if the vector $(\xi_0,\dots,\xi_k)^{\mathsf{T}}$ is an eigenvector of $-\eta H$ with eigenvalue $\lambda$, where $^{\mathsf{T}}$ denotes transpose.

When $\epsilon=0$ and $\eta=k$, we have $-\eta H=\operatorname{diag}(k,k-1,\dots,1,k)$, which has eigenvalue $k$ with multiplicity two.
The vector $(a,0,\dots,0,1)^{\mathsf{T}}$ is an eigenvector with eigenvalue $k$, and hence
$f_j(0,k,a,0,\dots,0,1,k)=0$ for $0\leqslant j\leqslant k+2$.
The Jacobian matrix of the functions $f_j$ $(0\leqslant j\leqslant k+2)$ with respect to the $k+3$ variables $\eta,\xi_0,\xi_1,\dots,\xi_k,\lambda$ is given by
\begin{equation*}
	\left(\begin{array}{c|cccc|c} p_0(\epsilon) \sum_{\ell=0}^k p_{\ell}(\epsilon) \xi_{\ell} &&&&& -\xi_0  \\ \vdots && \multicolumn{2}{c}{-\eta H-\lambda I} && \vdots \\ p_k(\epsilon) \sum_{\ell=0}^k p_{\ell}(\epsilon) \xi_{\ell} &&&&& -\xi_k \\ \hline 0 & 1 & 0 & \cdots & 0 & 0 \\ 0 & p_0(\epsilon) & p_1(\epsilon) & \cdots & p_k(\epsilon) & 0 \end{array}\right),
\end{equation*}
where $I$ denotes the identity matrix of degree $k+1$.
At $(0,k,a,0,\dots,0,1,k)$, this becomes
\begin{equation*}
	\left(\begin{array}{c|ccccc|c} 0 & 0 &&&&& -a  \\ 0 && -1 &&&& 0 \\ \vdots &&& \ddots &&& \vdots \\ 0 &&&& 1-k && 0 \\ 1 &&&&& 0 & -1 \\ \hline 0 & 1 & 0 & \cdots & 0 & 0 & 0 \\ 0 & 0 & 0 & \cdots & 0 & 1 & 0 \end{array}\right),	
\end{equation*}
which is nonsingular since $a\ne 0$, and hence the implicit function theorem applies:
there exist analytic functions $\eta(\epsilon),\xi_0(\epsilon),\dots,\xi_k(\epsilon),\lambda(\epsilon)$ of $\epsilon$ such that, on some neighborhood of $(0,k,a,0,\dots,0,1,k)$, we have $f_j(\epsilon,\eta,\xi_0,\xi_1,\dots,\xi_k,\lambda)=0$ for $0\leqslant j\leqslant k+2$ if and only if $\eta=\eta(\epsilon),\xi_0=\xi_0(\epsilon),\dots,\xi_k=\xi_k(\epsilon),\lambda=\lambda(\epsilon)$.
In fact, we have
\begin{equation}\label{xi_0}
	\xi_0(\epsilon)\equiv a.
\end{equation}
Moreover, it follows that
\begin{equation}\label{system of equations}
	\lambda(\epsilon)\xi_{\ell}(\epsilon)=r_{\ell}(\epsilon) \xi_{\ell}(\epsilon)+\eta(\epsilon)p_{\ell}(\epsilon),
\end{equation}
for $0\leqslant \ell\leqslant k$.
Observe that
\begin{equation}\label{evaluate p_l}
	p_{\ell}(\epsilon)=\sqrt{\frac{k!}{\ell!}}\,\epsilon^{k-\ell}+O(\epsilon^{k-\ell+1}).
\end{equation}
Set $\ell=0$ in \eqref{system of equations}.
Then since $\eta(0)=k$, it follows from \eqref{xi_0} and \eqref{evaluate p_l} that
\begin{equation}\label{lambda}
	\lambda(\epsilon)=r_0(\epsilon)+\frac{\sqrt{k!}\,k}{a}\epsilon^k+O(\epsilon^{k+1}).
\end{equation}
Note that
\begin{equation}\label{r0-rl}
	r_0(\epsilon)-r_{\ell}(\epsilon)=\ell+O(\epsilon^2),
\end{equation}
for $0\leqslant \ell\leqslant k$.
In particular, for $1\leqslant \ell\leqslant k$, $\epsilon=0$ is not a pole of $1/(\lambda(\epsilon)-r_{\ell}(\epsilon))$, and by \eqref{system of equations} we have
\begin{equation*}
	\xi_{\ell}(\epsilon)=\frac{\eta(\epsilon)p_{\ell}(\epsilon)}{\lambda(\epsilon)-r_{\ell}(\epsilon)}.
\end{equation*}
Hence
\begin{equation}\label{solve for eta}
	1=\sum_{\ell=0}^k p_{\ell}(\epsilon) \xi_{\ell}(\epsilon)=a p_0(\epsilon)+\eta(\epsilon)\sum_{\ell=1}^k\frac{p_{\ell}(\epsilon)^2}{\lambda(\epsilon)-r_{\ell}(\epsilon)}.
\end{equation}
Since
\begin{equation*}
	\frac{1}{\lambda(\epsilon)-r_{\ell}(\epsilon)}=\frac{1}{r_0(\epsilon)-r_{\ell}(\epsilon)}-\frac{\sqrt{k!}\,k}{\ell^2a}\epsilon^k +O(\epsilon^{k+1})
\end{equation*}
by \eqref{lambda} and \eqref{r0-rl}, it follows from \eqref{evaluate p_l} that
\begin{equation*}
	\sum_{\ell=1}^k\frac{p_{\ell}(\epsilon)^2}{\lambda(\epsilon)-r_{\ell}(\epsilon)}=\sum_{\ell=1}^k\frac{p_{\ell}(\epsilon)^2}{r_0(\epsilon)-r_{\ell}(\epsilon)}-\frac{\sqrt{k!}}{ka}\epsilon^k+O(\epsilon^{k+1}).
\end{equation*}
Since the constant term above equals $1/k$ by \eqref{evaluate p_l} and \eqref{r0-rl}, it follows from \eqref{solve for eta} that
\begin{align*}
	\eta(\epsilon) &=\left(\sum_{\ell=1}^k\frac{p_{\ell}(\epsilon)^2}{\lambda(\epsilon)-r_{\ell}(\epsilon)}\right)^{\!\!\!-1}\!(1-a p_0(\epsilon)) \\
	&= \left[\left(\sum_{\ell=1}^k\frac{p_{\ell}(\epsilon)^2}{r_0(\epsilon)-r_{\ell}(\epsilon)}\right)^{\!\!\!-1}\!\!+\frac{\sqrt{k!}\,k}{a}\epsilon^k+O(\epsilon^{k+1})\right] \!\!(1-a p_0(\epsilon)) \\
	&= \left(\sum_{\ell=1}^k\frac{p_{\ell}(\epsilon)^2}{r_0(\epsilon)-r_{\ell}(\epsilon)}\right)^{\!\!\!-1}\!\!+\sqrt{k!}\,k\left(\frac{1}{a}-a\right)\!\epsilon^k+O(\epsilon^{k+1}).
\end{align*}
We note that the terms of degree less than $k$ in $\eta(\epsilon)$ are independent of $a$.

\section{Spatial search algorithm}\label{sec:search}

We now move on to the discussions on the search algorithm.
Recall that the Hamiltonian $H=-\gamma A-|w\rangle\langle w|$ is determined by $\epsilon=1/\sqrt{n}$ and $\eta=1/\gamma n=\epsilon^2/\gamma$.
We consider three kinds of $\eta$ as follows.
In the preceding discussions, we choose $a=\pm 1$, so that
\begin{equation*}
	\frac{1}{a}-a=0.
\end{equation*}
We may choose any real scalar $\rho$ and two values of $a$ so as to $1/a-a=\rho$ to achieve a quadratic speedup, but it can be shown that setting $\rho=0$ gives asymptotically optimal success probability $1+o(1)$.
We will use $^+$ (resp.~$^-$) to mean functions associated with $a=+1$ (resp.~$a=-1$), e.g., $\lambda^+,\xi_i^+$.
In particular, we have $\eta^{\pm}$ and thus two Hamiltonians $H^{\pm}$.
For the search algorithm, we will use the Hamiltonian $H=H^{\circ}$ determined when $\eta$ is
\begin{equation}\label{choice of eta}
	\eta^{\circ}:=\left(\sum_{\ell=1}^k\frac{p_{\ell}(\epsilon)^2}{r_0(\epsilon)-r_{\ell}(\epsilon)}\right)^{\!\!\!-1}\!\!.
\end{equation}
Note that
\begin{equation}\label{three etas are close}
	\eta^{\pm}=\eta^{\circ}+O(\epsilon^{k+1}).
\end{equation}

Let $\xi^{\pm}=(\xi^{\pm}_0,\dots,\xi^{\pm}_k)^{\mathsf{T}}$, so that
\begin{equation}\label{two eigenvectors}
	-\eta^{\pm}H^{\pm} \xi^{\pm}=\lambda^{\pm}\xi^{\pm}.
\end{equation}
We have
\begin{equation*}
	\xi^{\pm}=(\pm 1,0,\dots,0,1)^{\mathsf{T}}+o(1),
\end{equation*}
where the Landau notation is used entrywise for vectors and matrices.
In particular,
\begin{equation}\label{almost orthonormal}
	\|\xi^{\pm}\|^2=2+o(1), \qquad \langle \xi^+|\xi^-\rangle=o(1).
\end{equation}

Recall that the algorithm begins in the state
\begin{equation}\label{approximate D}
	\psi(0)=|{{s}}\rangle=|\lambda_0\rangle=(1,0,\dots,0)^{\mathsf{T}}=\frac{1}{2}\xi^+-\frac{1}{2}\xi^-+o(1).
\end{equation}
On the other hand, it follows from \eqref{evaluate p_l} that
\begin{eqnarray}\label{approximate w}
	|w\rangle &=& \sum_{\ell=0}^kP_{\ell}|w\rangle=\sum_{\ell=0}^kp_{\ell}(\epsilon)|\lambda_{\ell}\rangle=(0,\dots,0,1)^{\mathsf{T}}+o(1)\nonumber \\
	&=&\frac{1}{2}\xi^++\frac{1}{2}\xi^-+o(1).
\end{eqnarray}

\subsection{Computational complexity}

In order to analyze the algorithm, we will need to estimate the differences of the operators $\mathrm{e}^{-\mathrm{i}H^{\pm}t}$ and $\mathrm{e}^{-\mathrm{i}H^{\circ}t}$.
To this end, we recall Wilcox's formula~\cite{Wilcox1967JMP} for the derivative of a matrix exponential function.
Let $M(\eta)$ be a square matrix function of a real variable $\eta$ which is of class $C^1$, and let $t$ be another real variable.
Then we have
\begin{equation}\label{Wilcox}
	\frac{\partial}{\partial\eta}\mathrm{e}^{M(\eta)t} = \int_0^t \mathrm{e}^{M(\eta)(t-s)} \frac{dM}{d\eta}(\eta) \,\mathrm{e}^{M(\eta)s}\,ds.
\end{equation}
This formula is proved, e.g., by observing that, for every fixed $\eta$, both sides above satisfy the linear differential equation
\begin{equation*}
	\frac{d\Phi}{dt}(t)=M(\eta)\Phi(t)+\frac{dM}{d\eta}(\eta)\, \mathrm{e}^{M(\eta)t}
\end{equation*}
with initial condition $\Phi(0)=0$.
We use Wilcox's formula \eqref{Wilcox} with $M(\epsilon,\eta)=-\mathrm{i}H(\epsilon,\eta)$ to obtain
\begin{equation*}
	\frac{\partial}{\partial\eta}\mathrm{e}^{-\mathrm{i}H(\epsilon,\eta)t} = -\mathrm{i}\int_0^t \mathrm{e}^{-\mathrm{i}H(\epsilon,\eta)(t-s)} \frac{\partial H}{\partial\eta}(\epsilon,\eta) \,\mathrm{e}^{-\mathrm{i} H(\epsilon,\eta)s}\,ds.
\end{equation*}
Since
\begin{equation*}
	H=-\frac{1}{\eta}\operatorname{diag}(r_0(\epsilon),r_1(\epsilon),\dots,r_k(\epsilon)) -\big(p_{\ell}(\epsilon)p_{\ell'}(\epsilon)\big)_{\ell,\ell'=0}^k,
\end{equation*}
we have
\begin{equation*}
	\frac{\partial H}{\partial\eta}(\epsilon,\eta)=\frac{1}{\eta^2} \operatorname{diag}(r_0(\epsilon),r_1(\epsilon),\dots,r_k(\epsilon)).
\end{equation*}
Assuming that $\eta$ is sufficiently close to $k$, the entries of this matrix are bounded, irrespective of $\epsilon=1/\sqrt{n}$.
The matrix $\mathrm{e}^{-\mathrm{i} H(\epsilon,\eta)s}$ is diagonalized by a real orthogonal matrix, and its eigenvalues all have modulus one, from which it follows that all of its entries are bounded by, say, $k+1$, irrespective of $\epsilon,\eta$, and $s$.
A similar result holds for $\mathrm{e}^{-\mathrm{i} H(\epsilon,\eta)(t-s)}$.
By these comments, all the entries of the integrand are uniformly bounded by a constant $R>0$, which is independent of the parameters (as long as $\eta$ is close to $k$).
It follows that the entries of $\partial \mathrm{e}^{-\mathrm{i}H(\epsilon,\eta)t}/\partial\eta$ are bounded by $R|t|$.
By the mean value theorem, it follows that the entries of $\mathrm{e}^{-\mathrm{i}H^{\pm}t}-\mathrm{e}^{-\mathrm{i}H^{\circ}t}$ are bounded by $R|t(\eta^{\pm}-\eta^{\circ})|$.

We set the running time as
\begin{equation*}
	t_{\mathrm{run}}=\frac{\pi n^{k/2}}{2\sqrt{k!}}=\frac{\pi}{2\sqrt{k!}\,\epsilon^k} \left(\approx \frac{\pi \sqrt{N}}{2} \right).
\end{equation*}
Since $t_{\mathrm{run}}=O(n^{k/2})=O(1/\epsilon^k)$, it follows from the above comments and \eqref{three etas are close} that
\begin{equation*}
	\mathrm{e}^{-\mathrm{i}H^{\pm}t_{\mathrm{run}}}=\mathrm{e}^{-\mathrm{i}H^{\circ}t_{\mathrm{run}}}+o(1).
\end{equation*}
We have also mentioned that $\mathrm{e}^{-\mathrm{i}H^{\pm}t_{\mathrm{run}}},\mathrm{e}^{-\mathrm{i}H^{\circ}t_{\mathrm{run}}}=O(1)$.
Hence it follows from \eqref{two eigenvectors} and \eqref{approximate D} that
\begin{align*}
	\psi(t_{\mathrm{run}}) &= \mathrm{e}^{-\mathrm{i}H^{\circ}t_{\mathrm{run}}}|{{s}}\rangle \\
	&= \frac{1}{2}\mathrm{e}^{-\mathrm{i}H^+t_{\mathrm{run}}}\xi^+ -\frac{1}{2}\mathrm{e}^{-\mathrm{i}H^-t_{\mathrm{run}}}\xi^- +o(1)\\
	&= \frac{\mathrm{e}^{\mathrm{i}\lambda^+t_{\mathrm{run}}/\eta^+}}{2}\xi^+ -\frac{\mathrm{e}^{\mathrm{i}\lambda^-t_{\mathrm{run}}/\eta^-}}{2}\xi^-+o(1).
\end{align*}
By \eqref{lambda} (applied to $a=\pm 1$), \eqref{three etas are close}, and since $\eta^{\circ}(0)=k$, we have
\begin{equation*}
	\frac{\lambda^{\pm}(\epsilon)}{\eta^{\pm}(\epsilon)}=\frac{r_0(\epsilon)}{\eta^{\circ}(\epsilon)} \pm\sqrt{k!}\,\epsilon^k +O(\epsilon^{k+1}),
\end{equation*}
so that
\begin{equation*}
	\left(\frac{\lambda^+}{\eta^+}-\frac{\lambda^-}{\eta^-}\right)\!t_{\mathrm{run}}=\pi+o(1).
\end{equation*}
Hence
\begin{equation*}
	\psi(t_{\mathrm{run}})=\mathrm{e}^{\mathrm{i}\lambda^+t_{\mathrm{run}}/\eta^+} \left(\frac{1}{2}\xi^++\frac{1}{2}\xi^-\right)+o(1),
\end{equation*}
and it follows from \eqref{almost orthonormal} and \eqref{approximate w} that the success probability is
\begin{equation*}
	p_{\mathrm{succ}}=|\langle w|\psi(t_{\mathrm{run}})\rangle|^2=1+o(1),
\end{equation*}
concluding the analysis of the computational complexity.

\subsection{Some final observations}
We remark that $\eta^{\circ}$ can in fact be replaced by any analytic function of $\epsilon=1/\sqrt{n}$ whose series expansion agrees with the right-hand side of \eqref{choice of eta} up to terms of degree $k$.
Our choice of the parameter $\gamma$ is
\begin{equation*}
	\gamma^{\circ}:=\frac{\epsilon^2}{\eta^{\circ}}=\epsilon^2\sum_{\ell=1}^k\frac{p_{\ell}(\epsilon)^2}{r_0(\epsilon)-r_{\ell}(\epsilon)},
\end{equation*}
but the terms of degree higher than $k+2$ can be discarded or changed arbitrarily.

We list some values of $\gamma^{\circ}$ for small $k$:
\begin{eqnarray*}
\gamma^{\circ}\Big\rvert_{k=3} &=& \dfrac{ \epsilon^2 (1-3\epsilon^2) ( 2 +\epsilon^2 +16
   \epsilon^4 -52 \epsilon^6 +24 \epsilon^8) }{ 6 (1-\epsilon^2)^2 (1-2\epsilon^2)^2 }, \\
\gamma^{\circ}\Big\rvert_{k=4} &=& \dfrac{ \epsilon^2 (1-4\epsilon^2) (  3 -11 \epsilon^2 +33 \epsilon^4 +47 \epsilon^6 -660 \epsilon^8 +1116 \epsilon^{10} -432 \epsilon^{12} ) }{ 12 (1-\epsilon^2)^2 (1-2\epsilon^2)^2 (1-3\epsilon^2)^2 }, \\
\gamma^{\circ}\Big\rvert_{k=5} &=&  \epsilon^2 (1-5\epsilon^2) \big(12 -117 \epsilon^2 +532 \epsilon^4 -1107 \epsilon^6 +2508 \epsilon^8 - \\
&& \hspace{2cm} 22588 \epsilon^{10}+80448\epsilon^{12} -99648 \epsilon^{14} +34560 \epsilon^{16}\big)/\\
&& \hspace{2cm}  { \left(60 (1-\epsilon^2)^2 (1-2\epsilon^2)^2 (1-3\epsilon^2)^2 (1-4\epsilon^2)^2\right). }
   \end{eqnarray*} 
It seems that there is no ``clean'' expression for $\gamma^{\circ}$.
For $k=3$, we have
\begin{equation*}
	\gamma^{\circ}=\frac{\epsilon^2}{3}+\frac{7\epsilon^4}{6}+O(\epsilon^6),
\end{equation*}
which was used by Wong~\cite{Wong2016JPA} in his algorithm for $J(n,3)$.

\section{Conclusions}\label{sec:conc}

We have analyzed the quantum spatial search algorithm on Johnson graphs $J(n,k)$ by continuous-time quantum walks. We have explored the distance-transitivity to define an invariant reduced subspace of the Hilbert space based on subsets of the vertex set that have the same distance to the marked vertex. The key point to have success in the calculation process is to use an alternative basis for the reduced space, which is obtained by projecting the marked state $|w\rangle$ onto the reduced space with projectors $P_\ell$, which are projectors onto the eigenspaces of the adjacency operator associated with eigenvalues $\lambda_\ell$, for $0\leqslant \ell \leqslant k$. Two eigenvectors of the Hamiltonian are used to obtain the time evolution of the quantum walk, and by taking a running time equal to $\pi\sqrt{N}/2$ approximately, we have showed that the success probability is $1+o(1)$ for every fixed $k$. Our proof is fully rigorous from the mathematical point of view and generalizes the continuous-time Grover search on the complete graph.

As a future work, we are interested in generalizing our method for larger subclasses of distance-transitive graphs or distance-regular graphs.

\section*{Acknowledgements}
The work of H. Tanaka was supported by JSPS KAKENHI grant number JP20K03551.
The work of R. Portugal was supported by FAPERJ grant number CNE E-26/202.872/2018, and CNPq grant number 308923/2019-7.


\end{document}